\documentclass[11pt]{amsart}
\usepackage{amssymb,amscd, graphicx, epsf}

\begin{document}

\newtheorem{defi}{Definition}
\newtheorem*{thm}{Theorem}
\newtheorem*{conjecture}{Conjecture}
\newtheorem{prop}{Proposition}
\newtheorem*{rem}{Remark}
\newtheorem{numrem}{Remark}
\newtheorem{lemma}{Lemma}
\newtheorem{coro}{Corollary}
\newtheorem{quest}{Question}
\newtheorem{corollary}{Corollary}

\title[Dense packing with various solids]{Dense packing of space with various convex solids}

\author{Andr\'{a}s Bezdek}
\address{A.~Bezdek, MTA R\'{e}nyi Institute, 13-15 R\'{e}altanoda u., Budapest, Hungary \hfill\break
\hbox{\hspace{.53in}} {\it also}: Mathematics \& Statistics, Auburn University, Auburn, AL 36849-5310, USA}
\email{bezdean@auburn.edu}

\thanks{Both authors gratefully acknowledge research support: A.~Bezdek was supported by the Hungarian Research Foundation OTKA, grant \#068398;  W.~Kuperberg was supported in part by the DiscConvGeo (Discrete and Convex Geometry) project, in the framework of the European Community's
``Structuring the European Research Area'' programme.}

\author{W{\l}odzimierz Kuperberg}
\address{W.~Kuperberg, Mathematics \& Statistics, Auburn University, Auburn, AL 36849-5310, USA}
\email{kuperwl@auburn.edu}

\dedicatory{Dedicated to Professor L\'{a}szl\'{o} Fejes T\'{o}th}


\begin{abstract}  One of the basic problems in discrete geometry is to determine the most efficient packing
of congruent replicas of a given convex set $K$ in the plane or in space. The most commonly used measure
of efficiency is density.  Several types of the problem arise depending on the type of isometries allowed for
the packing:  packing by translates, lattice packing, translates and point reflections,  or all isometries.  Due to
its connections with number theory, crystallography, etc., lattice packing has been studied most extensively.
In two dimensions the theory is fairly well developed, and there are several significant results on lattice
packing in three dimensions as well.  This article surveys the known results, focusing on the most recent progress.
Also, many new problems are stated, indicating directions in which future development of the general packing
theory in three dimensions seems feasible.  

\end{abstract}  

\maketitle

\section{Definitions and Preliminaries.}

A $d$-dimensional {\em convex body} is a compact convex subset of ${\mathbb R}^n$, contained in a $d$-dimensional
flat and with non-void interior relative to the flat.  A $2$-dimensional convex body is called a {\em convex disk}.
The ($d$-dimensional) volume of a $d$-dimensional convex body $K$ will be denoted by ${\rm Vol}(K)$, but for $d=2$
we will sometimes alternately use the term ``area'' and the notation ${\rm Area}(K)$.

The {\em Minkowski sum} of sets $A$ and $B$ in ${\mathbb R}^d$ is defined as the set $$A+B=\{x+y: x\in A,\,y\in B\}.$$
If $A$ consists of a single point $a$, we write simply $a + B$ instead of $\{ a \} + B$.

For every convex body $K$ in ${\mathbb R}^d$ and every real number $\lambda$, the set $\lambda K$ is defined as
$\{\lambda x: x\in K\}$.  We usually write $-K$ instead of $(-1)K$, and $K-L$ instead of $K+(-1)L$.  A convex body $K$
in ${\mathbb R}^d$ is {\em centrally symmetric} if there is a point $c\in{\mathbb R}^d$ (the {\em center} of $K$) such that
$K=2c-K$.  For each convex body $K$, the centrally symmetric convex body ${\bf D}K = \frac{1}{2}(K-K)$ is called {\em the difference
body} of $K$.  

A {\em packing} of ${\mathbb R}^d$ is a family of $d$-dimensional convex bodies $K_i$ whose interiors are mutually
disjoint. A packing is a {\em tiling} if the union of its members is the whole space ${\mathbb R}^d$.

In what follows, we consider mostly packings with congruent replicas of a convex body $K$.
If the family ${\mathcal P}=\{K_i\}$ $(i=1,2,\ldots)$ of congruent replicas $K_i$  of a $d$-dimensional convex body $K$
is a packing, then {\em density} of  ${\mathcal P}$ is defined as
$$
d({\mathcal P})  =\limsup_{r\to\infty} \frac{1}{{\rm Vol}({\bf B}(r))}\sum_{i=1}^\infty{\rm Vol}(K_i\cap{\bf B}(r)),
$$
where ${\bf B}(r)$ is the ball of radius $r$, centered at the origin.  The supremum of $d({\mathcal P})$ taken over all
packings ${\mathcal P}$ with congruent replicas of $K$ is called {\em the packing density of $K$} and is denoted
by $\delta(K)$. The supremum is actually the maximum, as a densest packing with replicas of $K$ exists (see Groemer \cite{G86}).  
In case the allowed replicas of $K$ are restricted to translates of $K$ or to translates of $K$ by a lattice, the corresponding
packing densities of $K$ are denoted by $\delta_T(K)$ and by $\delta_L(K)$, respectively.  We also consider packings in which
translates of $K$ and translates of $-K$ are used; the corresponding packing density is denoted by $\delta_{T^{\text{*}}}(K)$. The
lattice-like version requires that each packing consists of translates of a non-overlapping pair $K\cup(v-K)$ by the vectors
of a lattice; the corresponding density is denoted by  $\delta_{L^{\text{*}}} (K)$ (here both the lattice $L$ and the vector $v$ are chosen
so that the resulting packing is of maximum density).  Naturally, the more restrictions are imposed on the type of the allowed
packing arrangements, the smaller is the corresponding packing density, therefore
$$0<\delta_L(K)\le\delta_T(K)\le\delta_{T^{\text{*}}}((K)\le\delta(K)\le1$$
and
$$
0<\delta_{L^{\text{*}}}(K)\le\delta_{T^{\text{*}}}((K)\le\delta(K)\le1.
$$

Obviously, if space ${\mathbb R}^d$ can be tiled by congruent replicas of $K$, then $\delta(K)=1$.  The converse is less
obvious, but not very difficult to prove: If $\delta(K)=1$, then ${\mathbb R}^d$ can be tiled by congruent replicas of $K$.
Similarly, if $\delta_T(K)=1$, then $K$ can tile space by its translated replicas; and if $\delta_{T^{\text{*}}}(K)=1$, then space can
be tiled by translates of $K$ combined with translates of $-K$.

It is well-known that a family  ${\mathcal P}=\{K+v_i \}$ of  translates of a
convex body $K$ is a packing if and only if the family ${\mathcal P}'=\{{\bf D}K+v_i \}$ is a packing (see \cite{M04}, also \cite{EGH89}, \cite{L69}
and \cite{GL87}).  This implies immediately that
$$
\delta_T(K)=\frac{{\rm Vol}(K)}{{\rm Vol}({\bf D}K)}\,\delta_T({\bf D}K)\le\frac{{\rm Vol}(K)}{{\rm Vol}({\bf D}K)}\,, \eqno(1.1)
$$
which gives a meaningful ({\em i.e.}, smaller than $1$) upper bound in case $K$ is not centrally symmetric.  The analogous statement and
bound hold for the lattice packing density $\delta_L$.

For more details, definitions, and basic properties on these notions, see \cite{FK93}.  For an overview of lattices and lattice packings, see
\cite{EGH89}, \cite{L69} and \cite{GL87}.

\section{Introduction.}

In contrast to the well developed theory of packing in two dimensions, there are not many results about packing densities of
convex bodies in ${\mathbb R}^3$.  With few exceptions, most of such results simply provide the value of the packing
density $\delta_L(K)$ for a specific convex body $K$, usually obtained by means of a classical method described by
Minkowski \cite{M04}.   In the next section we review those results, occasionally citing and describing some relevant results
about packing the plane ${\mathbb R}^2$ with congruent replicas of a {\em convex disk} (a convex body of dimension $2$).

In Sections 6, 7, and 8 we consider two simple types of convex bodies in ${\mathbb R}^3$, namely cones and cylinders. 
Given a convex disk $K$ in ${\mathbb R}^3$ and a point $v$ not in the plane of $K$, the {\em cone with base $K$ and
apex $v$}, denoted by  ${\bf C}_v(K)$, is the union of all line segments with one end at $v$  and the other one in $K$.
Given a convex disk $K$ in ${\mathbb R}^3$ and a line segment $s$ not parallel to the lane of $K$, the {\em cylinder with
base $K$ and generating segment $s$}, denoted by ${\bf\Pi}_s(K)$ is the Minkowski sum $s+K$.  (Observe that, with the exception
of tetrahedra, the base and apex of a cone are uniquely determined by the cone itself; likewise, with the exception of the
parallelepipeds, a cylinder has two bases exactly - one is a translate of the other, and its generating segment is determined
uniquely up to translation.) These two simple types of convex bodies we suggest to investigate first as a first step towards
building a systematic theory of packing in dimension three.  The plan is particularly suitable for the study of densities
$\delta_T$, $\delta_L$, $\delta_{T^{\text{*}}}(K)$, and $\delta_{L^{\text{*}}}(K)$, because of the affine invariance of the corresponding
problems.  Both for the cone and for the cylinder, each of the packing densities mentioned above depends only on the affine class
of the base. In Section 5 we describe in detail the nature of the affine invariance, we draw some immediate conclusions concerning those
suitable densities, and we state a few fundamental  open problems.

\section{Lattice Packing in Space.}

We begin with the following table listing a few convex bodies in ${\mathbb R}^3$ whose lattice packing densities $\delta_L$ have been explicitly
computed.

\vspace{5mm}

\begin{center}

\begin{tabular}{|r| c | c | l |}  \hline
\rule[-4mm]{0mm}{10mm} \bf \# & \bf Body & \bf Packing Density $\delta_L$ & \bf Author \& Reference\\
\hline \hline
\rule[-4mm]{0mm}{9mm} 1 & Ball $\{ x: |x|\le1\}$ &  $\frac{\pi}{\sqrt{18}}=0.74048\ldots$ & Gauss \cite{G840}\\
\hline
\rule[-3mm]{0mm}{8mm} 2 & Regular octahedron & $\frac{18}{19}=0.9473\ldots $ & Minkowski \cite{M04}\\
\hline
\rule[-2mm]{0mm}{6mm}  &  &  & Chalk \& Rogers \cite{CR48}, \\
\rule[-2mm]{0mm}{6mm} \raisebox{1.5ex}[0pt] 3 & \raisebox{1.5ex}[0pt]{Cylinder $C={\bf\Pi}_s(K)$} & \raisebox{1.5ex}[0pt]{$\delta_L(C)=\delta_L(K)$} &
\raisebox{.5ex}[0pt]{also Yeh \cite{Y48}}\\
\hline
\rule[-2mm]{0mm}{6mm}   &  Slab of a cube & {} &{}\\
\rule[-2mm]{0mm}{6mm} \raisebox{1.5ex}[0pt] 4 & \raisebox{.3ex}[0pt]{(see definition below)} & \raisebox{1.5ex}[0pt] {(see formula below)} &
\raisebox{1.5ex}[0pt]{Whitworth \cite{Wh48}}\\
\hline
\rule[-2mm]{0mm}{6mm}  &  Slab of a ball &{} &{} \\
\rule[-2mm]{0mm}{6mm} \raisebox{1.5ex}[0pt] 5 & \raisebox{.3ex}[0pt]{(see definition below)} & \raisebox{1.5ex}[0pt] {(see formula below)} &
\raisebox{1.5ex}[0pt]{Chalk \cite{C50}}\\
\hline
\rule[-2mm]{0mm}{6mm}  & Double cone & {}  & {}  \\
\rule[-2mm]{0mm}{6mm} \raisebox{1.5ex}[0pt] 6 & \raisebox{.3ex}[0pt]{(see definition below)} & \raisebox{1.5ex}[0pt]{$ \pi\sqrt6/9=0.85503\ldots$} &  \raisebox{1.5ex}[0pt]{Whitworth \cite{Wh51}}\\
\hline
\rule[-3mm]{0mm}{8mm} 7 & Tetrahedron & $\frac{18}{49}=0.3673\ldots $ & Hoylman \cite{H76}\\
\hline

\end{tabular}

\vspace{4mm}

Table 1.

\vspace{2mm}


\begin{figure}[h]
\includegraphics[scale=.5]{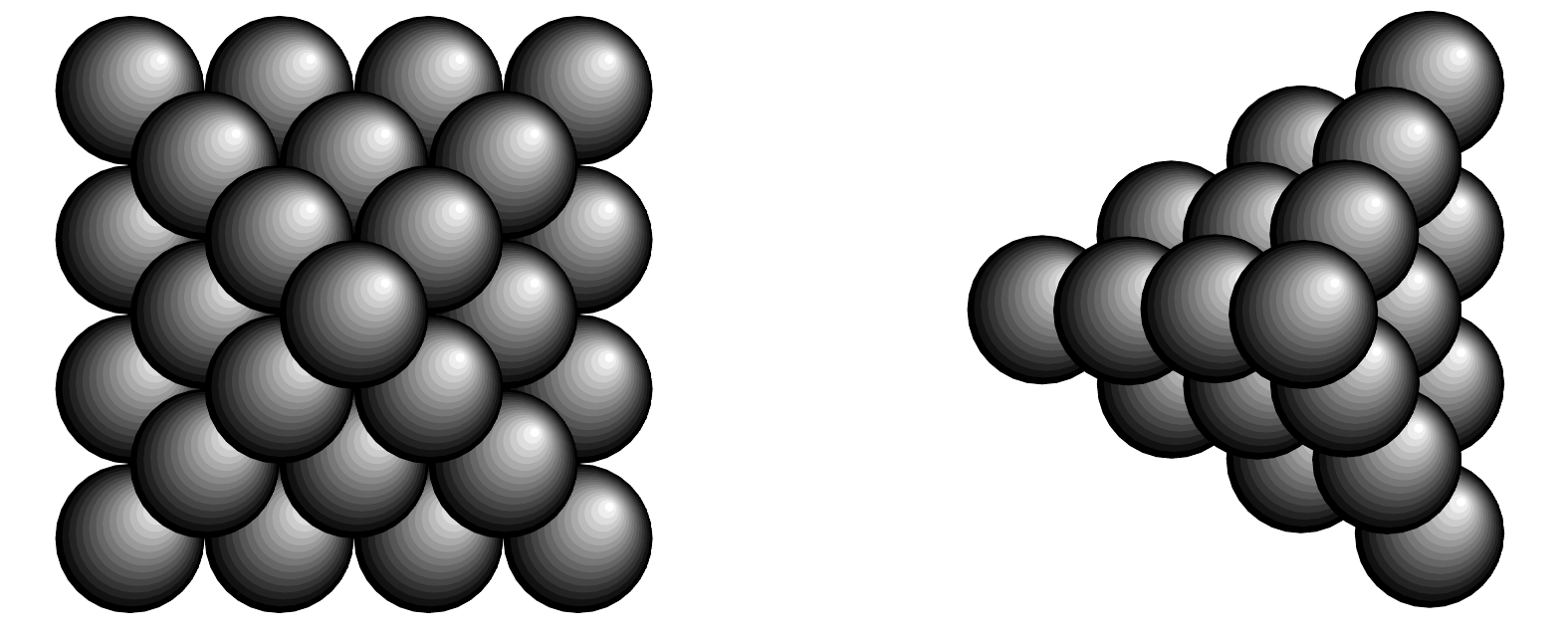}
\caption{Two clusters in the densest lattice packing of balls:  a ``square pyramid'' and a ``regular tetrahedron.'' }
\label{balls}
\end{figure}
\end{center}

{\bf Comments to Table 1.}

\begin{enumerate}

\item[1.] The densest lattice arrangements of spheres (balls) in ${\mathbb R}^3$ (see Fig.~\ref{balls}) was described already by Kepler \cite{K611},
but unsupported by proof, Kepler's assertion can only be considered a conjecture. The first one to prove that
$\delta_L({\bold B}^3)= \frac{\pi}{\sqrt{18}}$ was Gauss \cite{G840}.  Actually, Kepler asserted that the lattice arrangement shown in Fig.~\ref{balls}
is of maximum density among {\em all} sphere packings. This stronger conjecture, however, turned out to be extremely difficult to prove (see
Section 8, subsection 8.1). 

\item[2.] The regular octahedron is also called the regular $3$-dimensional {\em cross-polytope} and is denoted by ${\bold X}^3$.  Using his method
for computing lattice packing density of a centrally symmetric convex body, Minkowski \cite{M04} proved that $\delta_L({\bold X}^3)=\frac{18}{19}$. 
He applied the same method to the tetrahedron, but without success, for in the process he made a mistake in assuming that the difference body of the
regular tetrahedron is the regular octahedron (see Comment 7 below). 
\vspace{2mm}
\item[3.] The seemingly obvious equality $\delta_L(C)=\delta_L(K)$ is not trivial at all.  The trivial part is the inequality $\delta_L(C)\ge\delta_L(K)$,
obtained by stacking layers of cylinders erected over the densest lattice packing of the plane with translates of the base, but the opposite inequality
is quite nontrivial, since a cross-section of a lattice packing of the cylinders by a plane parallel to the cylinders' bases need not be a lattice packing
of the bases in the plane, and,  {\em a priori}, the density of such a packing could be greater than in any lattice packing.  A result of L.~Fejes T\'{o}th \cite{F50}, independently discovered also by Rogers \cite{R51},  says that this in fact cannot happen, {\em i.e.},  the density of a packing
with translates of a convex disk cannot exceed the maximum density attained in a lattice arrangement.
\vspace{2mm}
\item[4.] The $\lambda$-slab of a cube ($0<\lambda\le3$) is defined as
$$K_\lambda=\{x\in{\mathbb R}^3: |x_i|\le1, i =1,2,3; |x_1+x_2+x_3|\le\lambda\},$$
and its lattice packing density is given by the formula
$$
\delta_L(K_\lambda)= 
\begin{cases} \frac{1}{9}(9-\lambda^2) & \text{if $0<\lambda\le\frac{1}{2}$,}
\\
\frac{1}{4}\lambda(9-\lambda^2)/(-\lambda^3-3\lambda^2+24\lambda-1) &\text{if $\frac{1}{2}\le\lambda\le1$,}
\\
\frac{9}{8}(\lambda^3 -9\lambda^2+27\lambda-3)/\lambda(\lambda^2-9\lambda+27) &\text{if $1\le\lambda\le1$.}
\end{cases}
$$

Whitworth uses Minkowski's method, and his result generalizes the case of the regular octahedron ($\lambda=1$), item $2$ in the Table.

\vspace{2mm}
\item[5.] The $\lambda$-slab of a ball ($0<\lambda\le1$) is defined as
$$B_\lambda=\{x\in{\mathbb R}^3: |x|\le1, |x_3|\le\lambda\},$$
and its lattice packing density is given by the formula
$$\delta_L(B_\lambda)=\frac{\pi}{6}\sqrt{3-\lambda^2}.$$

Chalk uses Minkowski's method, and his result generalizes the case of the ball ($\lambda=1$), item 1 in the Table.

\vspace{2mm}
\item[6.] The double cone (see Fig.~\ref{dblecone}) is the set
$$K=\{x\in{\mathbb R}^3: \sqrt{x_1^2+x_2^2} + |x_3| \le 1\}.$$

As in item 4, Whitworth uses Minkowski's method to establish the lattice packing density of $K$.

\begin{center}
\begin{figure}[h]
\includegraphics[scale=.55]{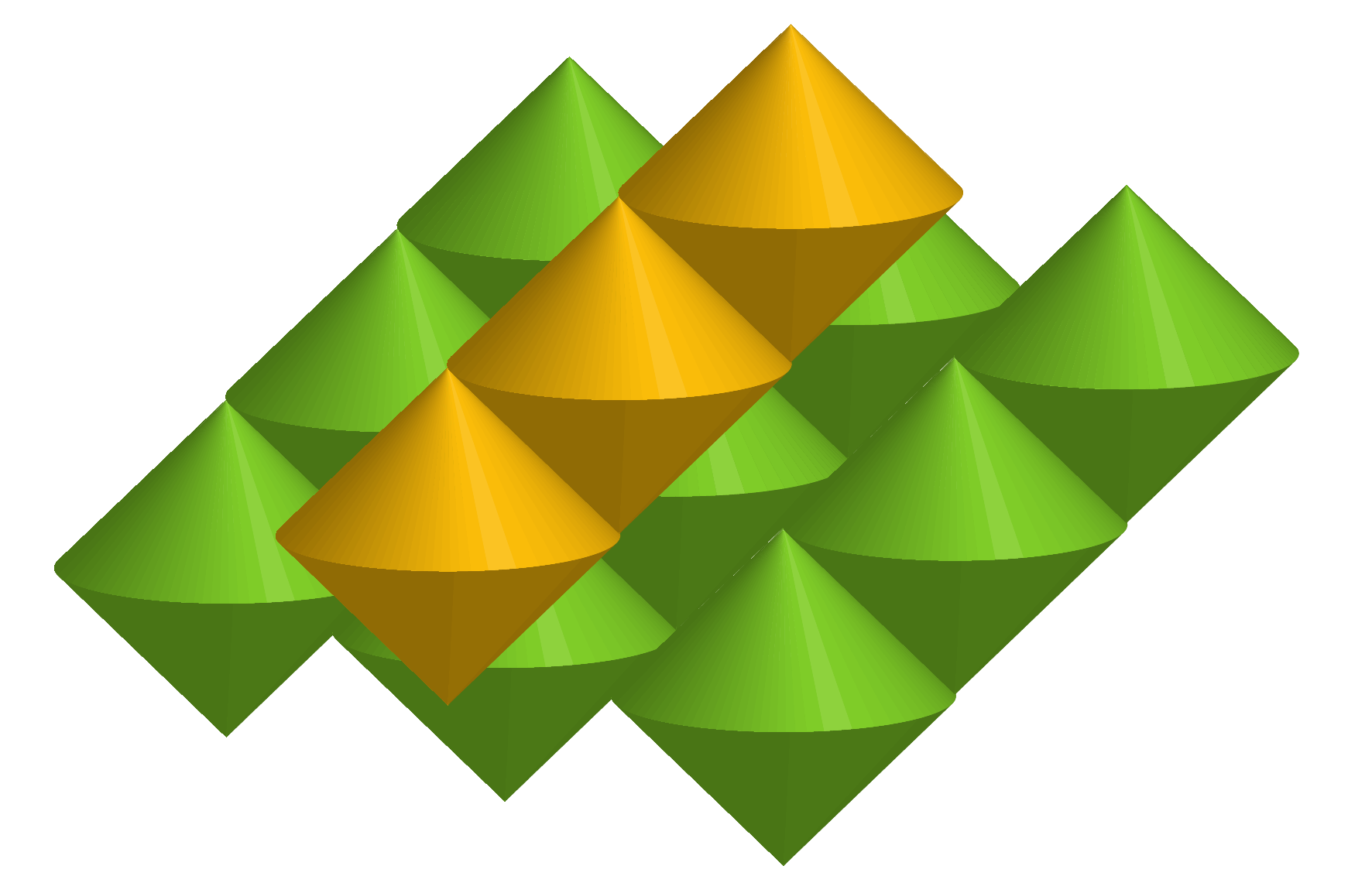}
\caption{The densest lattice packing with the double cone.}
\label{dblecone}
\end{figure}
\end{center}

\item[7.]   Minkowski's error in computing the lattice packing density of the tetrahedron was noticed by Groemer \cite{G62}, who proved
that $\frac{18}{49}$ is a lower bound for the density. Then  Douglas and Hoylman proved that Groemer's bound is in fact
the tetrahedron's lattice packing density.  The problem of the maximum density packing with congruent regular tetrahedra (allowing all
isometries) remains open and appears to be extremely difficult.  We report on the recent progress in the Section 8, subsection 8.4.

\begin{center}
\begin{figure}[h]
\includegraphics[scale=.85]{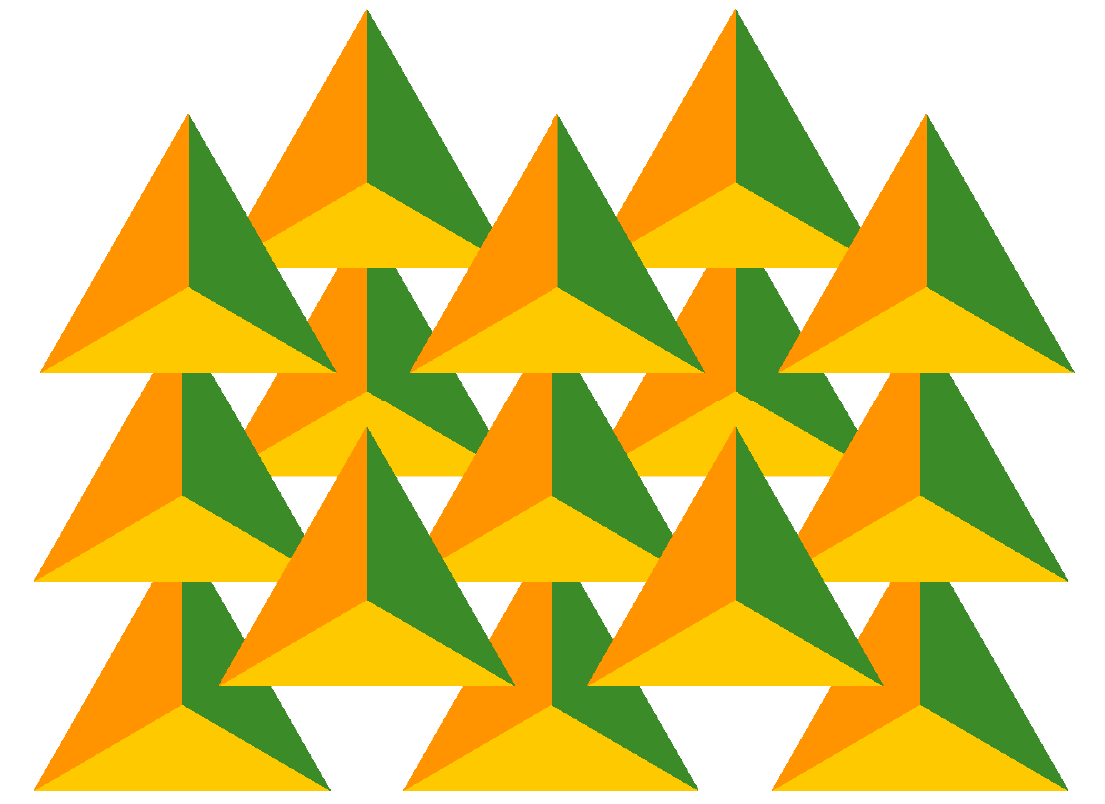}
\caption{The densest lattice packing with the tetrahedron.}
\label{tetra}
\end{figure}
\end{center}

\vspace{2mm}
\item[8.]  Each of the results listed in the table is obtained ``by hand,'' and, with the exception of Gauss, each of the authors uses Minkowski's
method.  The method often requires tedious computations with a large number of cases to analyze, which for some convex bodies becomes
prohibitively complex.  With the emergence of computer technology, however, it became possible to accomplish many such tasks in a very
short time.  In an impressive article published in year 2000, Betke and Henk \cite{BH00} present a fairly fast computer algorithm implementing
Minkowski's method for finding the lattice packing density of any $3$-dimensional convex polytope.  To show the algorithm's efficiency, the
article lists lattice packing density of each of the regular and Archimedean polytopes, many of which would be practically impossible to handle
without computers.

\end{enumerate}

\section{Packing Convex Bodies  by Translations.}

Thus far no example of a convex body $K$ has been found for which $\delta_T(K)>\delta_L(K)$.  In fact, there are only a few types of
convex bodies $K$ whose packing density $\delta_T(K)$ is known, namely:

\begin{enumerate}
\vspace{2mm}
\item any convex polytope $P$ that admits a tiling of space by its translates (it is known that each such polytope tiles space in a
lattice-like manner, in every dimension, see Venkov \cite{V54} or McMullen \cite{M80});
\vspace{2mm}
\item  any cylinder ${\bf C}_s K$ with a convex base $K$, since obviously $\delta_T({\bf C}_s K)= \delta_T(K)$;
\vspace{2mm}
\item  any non-symmetric body $K$ for which the packing density of the difference body $\delta_T({\bf D} K)$ is known.  For example,
the difference body of a body $K$ of constant width is a ball, hence the packing density of the ball can be used to find $\delta_T(K)$;
\vspace{2mm}
\item  any convex body $K$ such that ${\bf B}^3\subset K\subset RhD$, where $RhD$ denotes the rhombic dodecahedron circumscribing
the unit ball ${\bf B}^3$, which is the Voronoi polytope associated with the densest lattice packing of ${\bf B}^3$.
\vspace{2mm}
\end{enumerate}

The last two items are based on Hales' confirmation of the Kepler Conjecture, stating that
$\delta({\bf B}^3)=\delta_L({\bf B}^3)={\rm Vol}({\bf B}^3)/{\rm Vol}(RhD)$.

\vspace{2mm}
In contrast, in ${\mathbb R}^2$ it is known that
$$
\delta(K)=\delta_L(K)\quad \text{for every centrally symmetric convex disk $K$}, \eqno{(4.1)}
$$
 which implies that
$$
\delta_T(K)=\delta_L(K)\ \ \text{for every convex disk \(K\)}, \eqno{(4.2)}
$$
see L.~Fejes T\'{o}th \cite{F50}.

While equation $(4.2)$ perhaps holds true for $3$-dimensional convex bodies as well, equation $(4.1)$ does not, as the following example
shows.

Let $P$ be the (slightly irregular) affine-regular octahedron in ${\mathbb R}^3$ with vertices of the form $(\pm1,\pm1,0)$ and $(0,0,\pm1)$.
It is easy to see that $P$ cannot tile space by translates alone, hence $\delta_T(P)<1$.  On the other hand, $P$ can tile the plane with translates of itself combined with translates of its copies rotated by $90^\circ$ about the coordinate axes. Therefore $\delta(P)=1$.

\begin{center}
\begin{figure}[h]
\includegraphics[scale=.75]{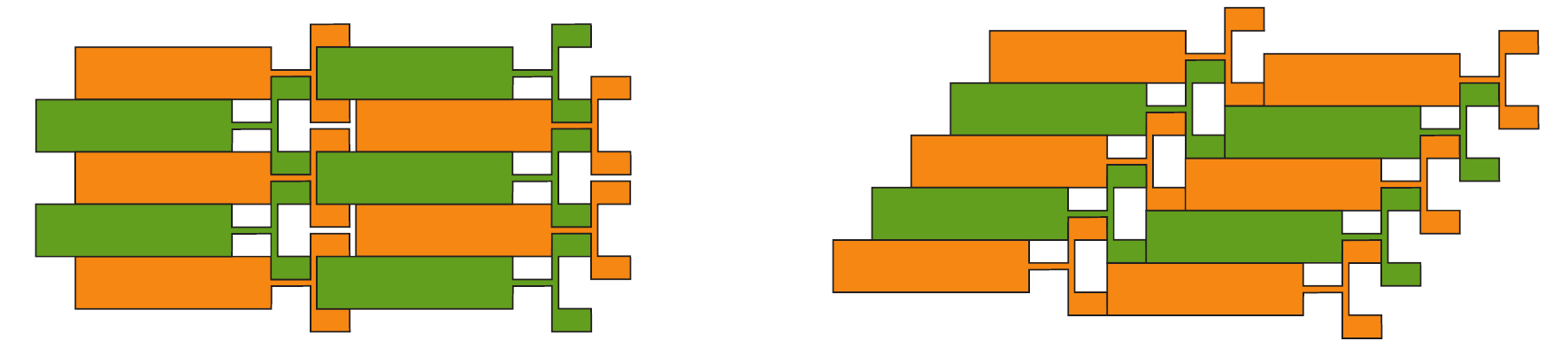}
\caption{An example of Bezdek and Kert\'{e}sz: a polygon whose translates can be packed more densely (left) than in its densest lattice packing (right).}
\label{poly}
\end{figure}
\end{center}

It should also be mentioned that already in dimension $2$ the assumption of convexity is indispensable for equation $(4.2)$.  A.~Bezdek and Kert\'{e}sz
 \cite{BK87} constructed a non-convex polygon that allows a dense non-lattice packing of the plane by its translates, denser than any lattice packing,
see Fig.~\ref{poly}.  (The construction of Bezdek and Kert\'{e}sz was modified by Heppes \cite{H90} so as to obtain a starlike polygon with the same property.)


The main question of this section remains open:

{\em Is it true that the maximum density of a packing with translates of a convex body in ${\mathbb R}^3$ is attained in a lattice packing?}

Similarly, the problem of whether or not $\delta_T^{\rm *}(K) = \delta_L^{\rm *}(K) $ holds for every $3$-dimensional convex body $K$ remains
open.

\section{Affine Invariance and Compactness.}

If $f: {\mathbb R}^d\to{\mathbb R}^d$ is an affine transformation, and if $K_1$ is a translate of a convex body $K$, then $f(K_1)$ is a translate
of $f(K)$.  Similarly, if $K_1$ is a translate of $-K$, then $f(K_1)$ is a translate of $-f(K)$.   Therefore the affine image of a packing with translates
of $K$ is a packing with translates of the image of $K$, and these two packings have the same density.  Moreover, the affine image of a lattice
packing with copies of $K$ is a lattice packing with the affine image of $K$.  The same {\em affine invariance} holds true for any packing that
combines translates of $K$ and of $-K$.   These simple facts imply that if convex bodies $K$ and $M$ are affine-equivalent, then:
$$
\delta_T(K)=\delta_T(M), \ \ \delta_{T^{\rm*}}(K)=\delta_{T^{\rm*}}(M), \ \ \delta_L(K)=\delta_L(M),\ \ \text{and\ \  $\delta_{L^{\rm*}}(K)=\delta_{L^{\rm*}}(M)$}.
$$

Therefore we can say that the domain of each of the four density functions $\delta_T$, $\delta_{T^{\rm*}}$, $\delta_L$, and  $\delta_{L^{\rm*}}$  is the
set of affine equivalence classes of convex bodies.  Let $[K]$ denote the affine equivalence class of the convex body $K$.  Following Macbeath
\cite{M51}, we supply the set of affine equivalence classes of convex bodies in
${\mathbb R}^d$ with the distance function $d$ defined as follows: for every pair $K, M$ of convex bodies, set 
$$
\rho(K,M)={\rm inf}\{{\rm Vol}(K')/{\rm Vol}(M): K' \text{ is affine equivalent to $K$ and } K'\supset M\}.
$$
Since the function $\rho$ is affine invariant,  the function $d$ given by
$$
d([K],[M])=\log\rho(K,M) + \log\rho(M,K)
$$
is well-defined. It is easy to check that $d$ is a metric on the set of all affine equivalence classes of convex bodies.  The space of such classes supplied
with this metric, denoted by ${\mathcal K}_a^d$, is compact (see Macbeath \cite{M51}), and each of the four packing density functions $\delta_T$,
$\delta_{T^{\rm*}}$, $\delta_L$, and  $\delta_{L^{\rm*}}$ defined on ${\mathcal K}_a^d$ is continuous.  Therefore each of them reaches its extreme
values.  Of course, the maximum value for each of them is $1$, reached at any convex body that tiles ${\mathbb R}^d$ by its translates.  However,
none of the four minimum values is presently known.

Determining those minimum values and the convex bodies at which they are attained seems
to be a very challenging problem, perhaps too difficult to expect to be solved in foreseeable future.  Reasonably good estimates for these minimum
values, however, should not be too hard to establish.

As for the maximum value of $1$, attained at the corresponding space tiling bodies (polytopes),
those that tile space by translations have been described in fairly simple terms by Venkov \cite{V54} and, independently, by McMullen \cite{M80}.
However, the analogous question, concerning which convex polytopes can tile space by their translates combined with translates of their negatives,
still remains unanswered.

\section{Packing Translates of Cones.}

We now turn our attention to the subspace $\mathcal C_a$ of ${\mathcal K}_a^3$ consisting of affine equivalence classes of cones, that is, affine classes
of bodies of the form  ${\bf C}_v(K)$, where $K$, the base, is a convex disk.  Since $\mathcal C_a$ is a closed subset of ${\mathcal K}_a^3$, it is
compact as well.  Notice that the affine class of the cone ${\bf C}_v(K)$ is determined uniquely by the affine equivalence class of its base $K$.
The affine class of a cone with base $K$ will be denoted by ${\bf C}K$.  Thus ${\bf C}K={\bf C}M$ if and ony if $K$ and $M$ are affinely equivalent convex
disks.

Again, the problem of maximum and minimum values arises that each of the four packing density functions attains on the compact set $\mathcal C_a$. 
This time, however, the maximum value of each of them is strictly smaller than $1$, since a cone cannot tile space, neither by its translates, nor by its translates combined with translates of its negative.  Thus we face a set of eight questions:

{\em Which convex disks produce cones of maximum and minimum packing density with respect to the four affine-invariant packing density functions?}

The eight extremum density values over the set of cones will be denoted by $c^{\rm max}$ and $c^{\rm min}$ supplied with the corresponding
subscripts $T$, $T^{\rm *}$,  $L$, and  $L^{\rm *}$.   The case of cones with centrally symmetric bases is of special interest, raising another set of eight analogous questions.

\begin{center}
\begin{figure}[h]
\includegraphics[scale=.75]{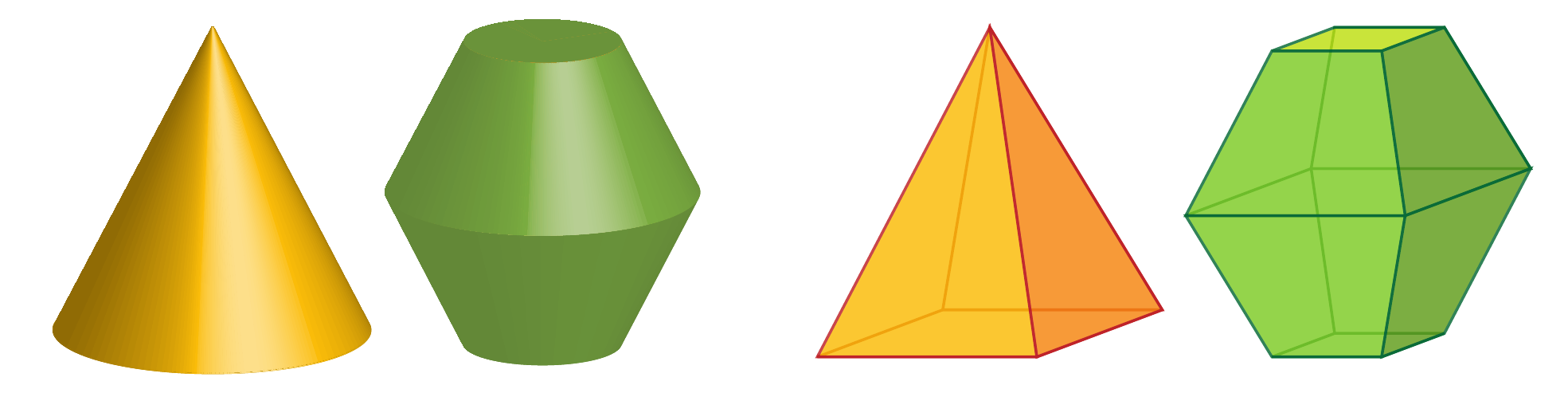}
\caption{Cones and their difference bodies:  the circular cone and the square pyramid.}
\label{diff}
\end{figure}
\end{center}

We begin with a lower bound for the volume of the difference body of a cone, to be used in the inequality $(1.1)$, producing an upper bound for the packing
density $\delta_T$ for all cones.  Figure \ref{diff} shows side-by-side two cones and their corresponding difference bodies.


For a cone with a centrally symmetric base, the volume ratio of the cone to its difference body is always $\frac{4}{7}$, which is easy to see.  For a cone
with non-symmetric  base, the corresponding volume ratio is always smaller than $\frac{4}{7}$, which follows directly from the Brunn-Minkowski inequality (see {\em e.g.} \cite{S93}) in dimension $2$.  The minimum ratio $\frac{2}{5}$ occurs for the triangular cone (the tetrahedron) only.  Thus we have the following upper bound:
$$
\delta_T({\bf C}K)<\frac{4}{7}
$$
as equality cannot occur since the difference body of any cone cannot tile space by translations.  Therefore
$$
c_T^{\rm max}=\max\{\delta_T({\bf C}K): K \text{ is a convex disk} \} < \frac{4}{7}.
$$

\begin{center}
\begin{figure}[h]
\includegraphics[scale=.8]{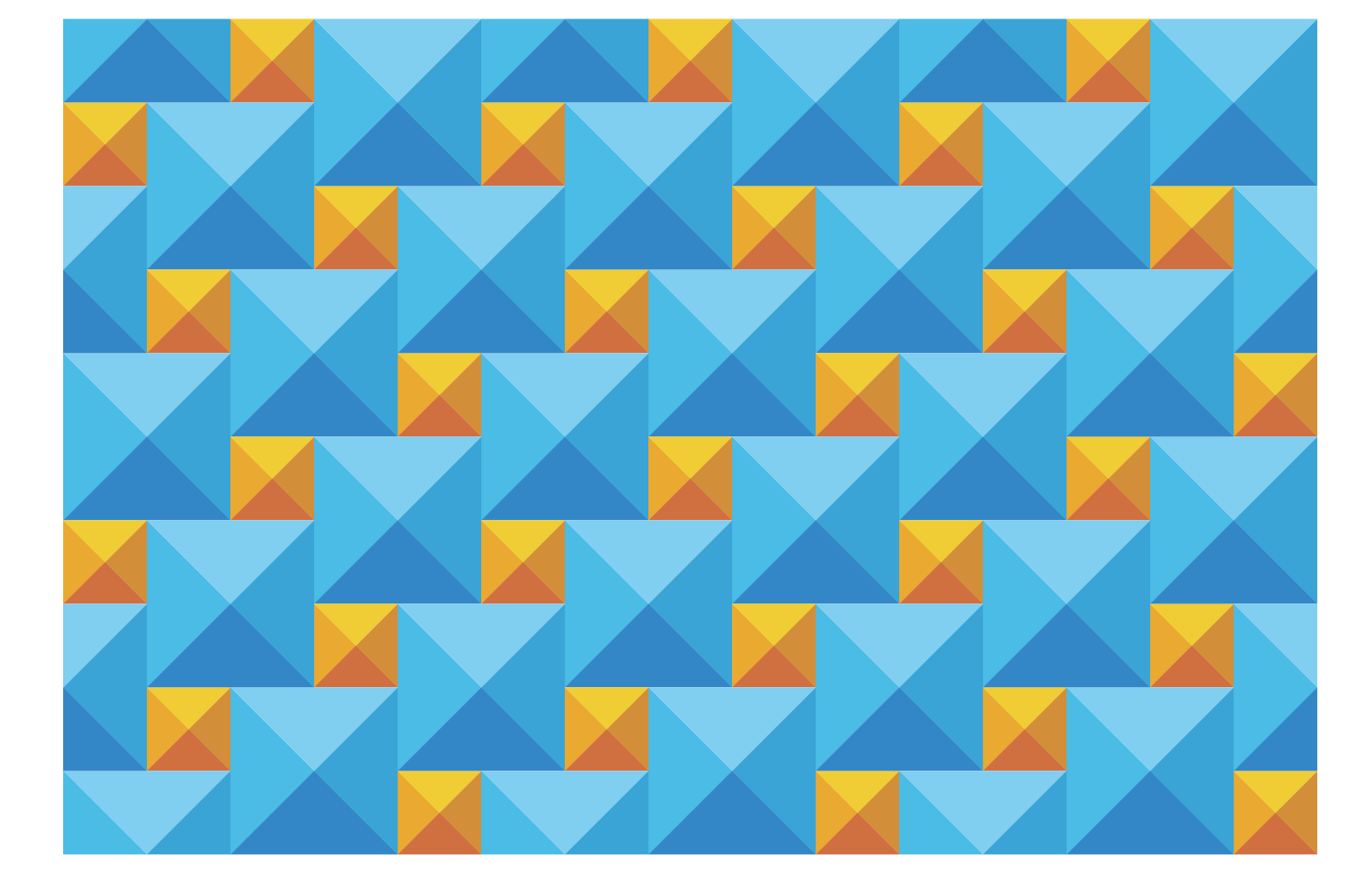}
\caption{A dense, though not the densest, lattice packing with the square pyramid.}
\label{square}
\end{figure}
\end{center}

On the other hand, there is a lattice packing with translates of a square pyramid, of density $\frac{8}{15}$.  The packing can be
described as follows.  Begin with a horizontal plane tiled by a lattice of ``L''-shaped figures consisting of a unit square with a
$\frac{1}{2}\times\frac{1}{2}$square attached to it.  Erect a square pyramid over each of the unit squares, get a layer of square pyramids,
in which the small squares are vacant. Place upon the first layer its translate, shifted so that the peaks of the pyramids form the first layer
plug the square holes in the second layer.  The vertical shift from the first layer to the second one is equal  to one-half of the pyramids'
height.  The two layers determine the entire lattice packing (see Fig.~\ref{square} for a top view of the two layers).


Thus $\delta_L({\bf C}S)\ge\frac{8}{15}$, where $S$ denotes the square.  However, according to the information supplied in private communication by
Betke and Henk, the authors of \cite{BH00}, the lattice packing density of the difference body of the square pyramid is
$\frac{112}{117}$, therefore $\delta_L({\bf C}S)\ge\frac{448}{819}=0.547\ldots >\frac{8}{15}=0.533\ldots$, and  we get the bounds

$$
0.547\ldots=\frac{448}{819}\le c_T^{\rm max}<\frac{4}{7}=0.571\ldots\,.
$$

\vspace{2mm}
{\bf Remark 1.}  By request of the authors of the present article,  Betke and Henk also computed the lattice packing densities
of the difference bodies of the pyramids with a regular hexagonal and a regular octagonal bases.  The results show that the lattice
packing density of the square pyramid is greater than those of the other two.  This seems to indicate that among all cones with centrally symmetric bases, the square-based cone has maximum lattice packing density.

\begin{center}
\begin{figure}[h]
\includegraphics[scale=1.25]{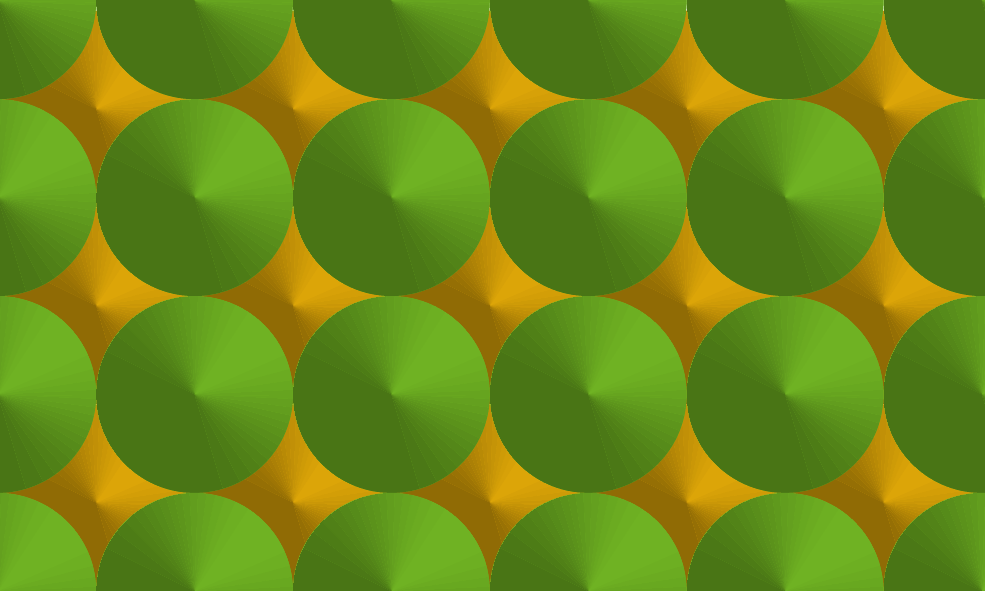}
\caption{A dense, though not likely the densest, lattice packing with the circular cone.}
\label{circ}
\end{figure}
\end{center}

{\bf Remark 2.}  The lattice packing density of the cone ${\bf C}E$ with a circular (elliptical) base $E$ has not been computed yet.  The best we know is
$$
 0.4469\ldots =\frac{2+\sqrt2}{24}\,\pi\le \delta_L({\bf C}E) \le \frac{\sqrt2}{9}\pi = 0.4936\ldots\,.
$$
The upper bound is found by inscribing a maximum volume ellipsoid in the difference body of the circular cone (see Fig.~\ref{diff}) and using the lattice packing density of the ball. The lower bound is obtained by the construction shown in Fig.~\ref{circ}.  Observe that the pattern is somewhat similar to that of the square pyramid seen in the previous figure.  Neither of the two bounds seems best possible - improvements should not be hard to obtain.


Besides the tetrahedron $T$, we do not know of any examples of cones with non-symmetric bases whose lattice packing density has
been computed.  The lattice packing density of the tetrahedron is  $\frac{18}{49}=0.3673\ldots $ (see Section 3, Table 1), which is perhaps the value of $c_L^{\rm min}$. 

\begin{center}

\begin{figure}[h]
\includegraphics[scale=.75]{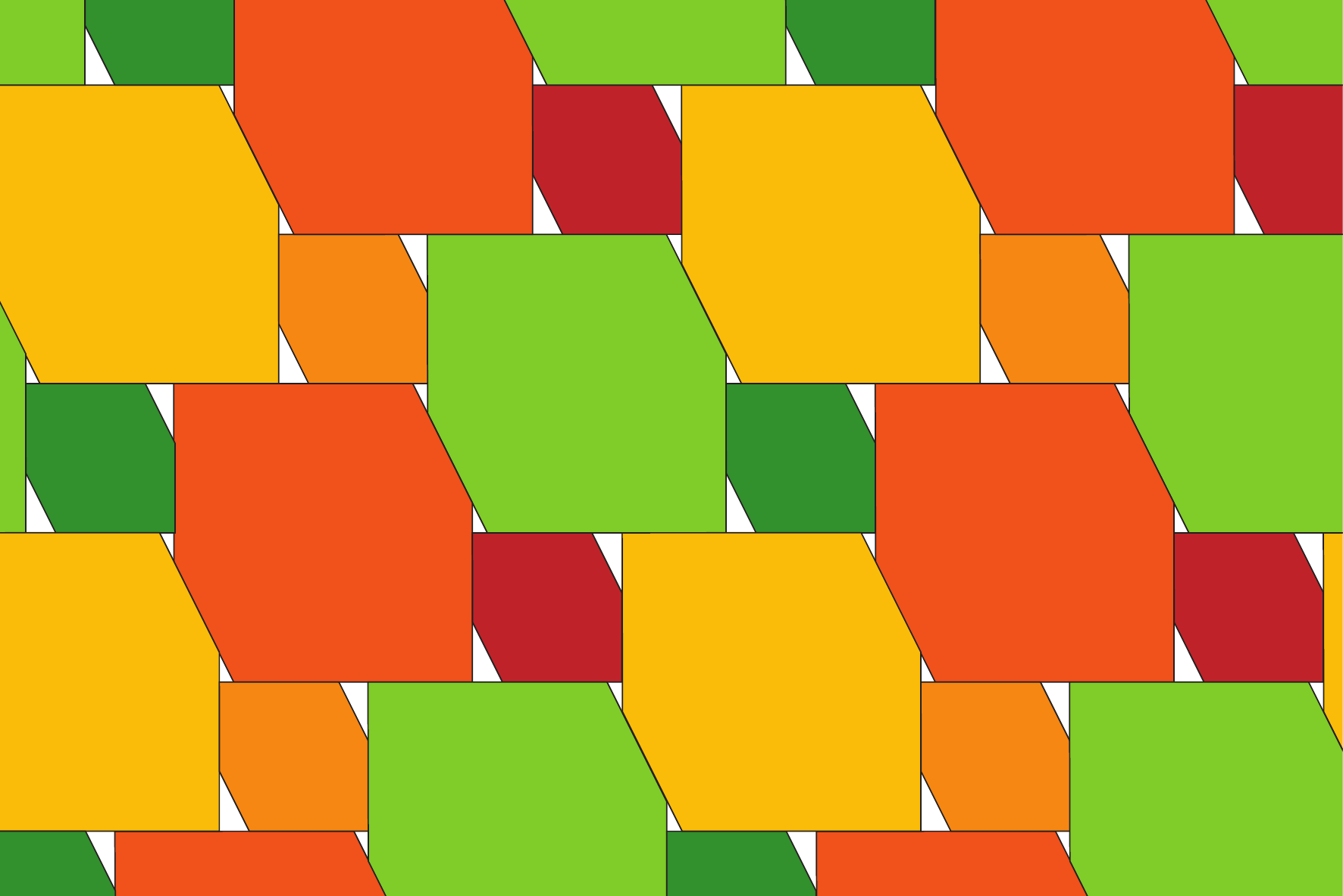}
\caption{A lattice packing with a pair of centrally symmetric hexagons $H$, $\frac{1}{2}H$.  The large hexagons form a packingof density at least
$\frac{3}{4}$.}
\label{hex} \end{figure}

\end{center}

Turning to cones with centrally symmetric bases, we obtain a common lower bound for their lattice packing density by a construction similar to that described for the square pyramid (see Fig.~\ref{square}).  First, observe that every centrally symmetric hexagon $H$ is contained in a parallelogram whose sides are extensions of sides of $H$ and of area at most $\frac{4}{3}$, maximum being reached by the regular hexagon. By an affine transformation we can assume that the parallelogram is a unit square, and $H$ is obtained by cutting off two congruent right triangles at two of its opposite corners.  Since the square is of minimum area among parallelograms containing $H$, the legs of the cut-off triangle cannot be longer than $\frac{1}{2}$. Therefore the arrangement shown in Fig.~\ref{hex} is a lattice packing of the plane with the pair consisting of $H$ and a translate of $\frac{1}{2}H$ attached to $H$, and the density of the collection of translates of $H$ (the large hexagons) is at least $\frac{3}{4}$, minimum being reached when $H$ is an affine regular hexagon. 

In a similar way as in the construction for the square pyramid, treating the small hexagons as holes in one layer of hexagonal pyramids, this packing gives rise to a lattice packing of space with the cone ${\bf C}H$. The density of this packing is at least $\frac{1}{2}$. 


Finally, by a theorem of Tammela \cite{T70}, every centrally symmetric convex disk $K$ of area $1$ is contained in a centrally 
symmetric hexagon $H$ of area at most $(3.570624)/4$, therefore, by the construction described above, we get the bound
$\delta_L({\bf C}K)\ge0.446328\ldots$ for every cone with a centrally symmetric disk $K$. Therefore 
$$
c_T^\text{min}\ge0.446328\ldots\,.
$$

\section{Packing Translates of Cones and Their Negatives.}

While we know that the value of  $c_{T^{\rm *}}^{\rm max}$ is smaller than $1$, an explicit upper bound below $1$ is not easy to produce. 
B\'{a}r\'{a}ny and Matou\v{s}ek \cite{BM07} found an explicit constant $\varepsilon>0$ such that the density of every packing of space with translates
of a cone and of its negative cannot exceed $1-\varepsilon$.  The value of $\varepsilon $ produced by their proof is very small, about
$10^{-42}$, and there seems to be room for improvement.

The ``best known'' case is the densest lattice packing of regular octahedra, of
density $\frac{18}{19}$ (see Table 1), showing that the constant $\varepsilon$ cannot be greater than
$\frac{1}{19}$, that is, $c_{L^{\rm *}}^{\rm max}\ge\frac{18}{19}$, but $c_{L^{\rm *}}^{\rm max}$ is very likely to be considerably greater than
$\frac{18}{19}$. Namely, it is likely that in the densest packing with translates of a square pyramid combined with translates of its negative, 
the pyramids do not form pairs joined by their common base.

Perhaps it is true in general that the maximum density of a packing with translates of a double cone with
a given centrally symmetric base is always smaller than some packing with translates of the cone and its negative. In other words, it seems likely
that separating the two parts of the double cone from each other always allows them to reach higher density.  It would be interesting to know
at least whether or not it is so for the double cone with a circular base and for the double square pyramid.

Elaborating on the idea of packing translates of the cone and its negative in pairs joined by their common, centrally symmetric, base, we use
a theorem of Petty \cite{P55}, stating that every centrally symmetric disk of area $1$ is contained in a parallelogram of area at most
$\frac{4}{3}$, the bound being sharp only in case of an affine regular hexagon.  This allows for enclosing such a pair of cones in an
affine regular octahedron whose densest lattice packing produces a packing with a cone with an arbitrary centrally symmetric base. The density
of so obtained packing is at least $\frac{3}{4}\times\frac{18}{19}=0.7105\ldots\,$.  A similar approach for cones with any convex base (the bases
of the cone and its negative need not coincide) produces a much weaker lower bound of $\frac{1}{2}\times\frac{18}{19}=0.47368\ldots$. (Here
the factor of $\frac{1}{2}$ is reached only in the case of the triangular base, that is, when the cone is a tetrahedron.) Thus we have
$$
\delta_{T^{\rm *}}({\bf C}K)\ge\frac{27}{38}=0.7105\ldots  \text{ for all cones with centrally symmetric convex bases $K$},  \eqno{(7.1)}
$$
and
$$
\delta_{T^{\rm *}}({\bf C}K)\ge\frac{9}{19}=0.47368\ldots \text{ for all cones with convex bases $K$}.  \eqno{(7.2)}
$$
The bound in $(7.1)$ is unlikely to be best possible, and the bound in $(7.2)$ definitely is not, since the construction of the presently known
densest packing with translates of the tetrahedron $\Delta$ and of $-\Delta$, recently found by Kallus,
Elser and Gravel \cite{KEG10}, is of density $(139+40\sqrt{10})/369=0.7194880\ldots\,$. 

\section{Packing Congruent Replicas of a Convex Body.} Here we consider packing densities $\delta(K)$ of a convex body $K$ in
${\mathbb R}^3$, with no restrictions on the nature of isometries used in packing.  There are not many bodies $K$ whose packing density
is known.  No good lower bound has been established for the packing density $\delta(K)$ valid for all convex $3$-dimensional bodies $K$.
A rather insignificant lower bound of $\frac{\sqrt3}{6}=0.288\ldots$ is easy to prove based on the known result in the plane, namely that the
packing density of every convex disk is at least $\frac{\sqrt3}{2}$  (see \cite{KK90}).

For centrally symmetric convex bodies $K$, the best known bound of this type is due to E.H.~Smith \cite{S05}, who proved that
$\delta(K)\ge0.53835\ldots$  body in ${\mathbb R}^3$.  The author indicates that the bound is not likely to be the best possible.  No reasonable
conjecture has been proposed, neither in the general case, nor under assumption of central symmetry, to  point to a specific convex body
whose packing density should be smallest among all convex bodies.  It is not even certain that such body exists.

Except for the trivial case of space-tiling polytopes, there are not many convex solids whose packing density $\delta$ (allowing all isometries)
is known.  In the following subsections we discuss known results for certain special cases.

\subsection{The Kepler conjecture.}  The three-dimensional sphere packing problem in its general form, without restrictions on the structure
of the spheres' arrangements is simple to state and easy to understand even for a non-expert.  The conjecture states that the maximum
density of a packing of  ${\mathbb R}^3$ with congruent balls is $\frac{\pi}{\sqrt{18}}=0.740480\ldots\,$, attained in the familiar lattice
arrangement (see Fig.~\ref{balls}).  The conjecture sounds very convincing to anyone who has ever seen  spherical objects, such as oranges
or apples, stacked in a pyramid, yet the proof eluded mathematicians for centuries. A problem so appealing attracts attention of experts and
laymen alike, and a solution tends to instantly elevate its author to the status of celebrity.  The Kepler conjecture, also known as the sphere
packing conjecture, has a long and fascinating history, see \cite{H00}.  The unsuccessful attempts at proof and the nature of the proof that was
produced at last seem to indicate that this is one of those problems that cannot be resolved with a reasonably simple and reasonably short proof.

The proof is due to Thomas Hales, who announced it in \cite{H98}, and then, during the past 13 years presented a series of articles on the subject
(see \cite{H97a, H97b, H05, H06a, H06b, H06c}, see also \cite{F06} by Ferguson, a student of Hales).   The description
of the theoretical approach to the problem and results of the work of computer occupies nearly 300 pages in these articles.  At the computational
stage of the proof, computers examined some $5,000$ computer-generated cases, each of the cases requiring optimization analysis of a system
of non-linear inequalities with a large number of variables. Hales main approach follows a strategy suggested by  L.~Fejes T\'{o}th in 1953
(see \cite{F72}) who anticipated a then insurmountable amount of computations needed for the case analysis.

As a corollary to Hales' result, the packing density of any convex body  $K$ such that ${\bf B}^3\subset K\subset RhD$ is easily computed:
$\delta(K)=\delta_T(K)=\frac{{\rm Vol}(K)}{{\rm Vol}(RhD)}$, where $RhD$ denotes the rhombic dodecahedron circumscribed about the unit ball
${\bf B}^3$.

In 2003, Hales launched a project called FLYSPECK, designed for an automatic (computerized) formal verification of his proof.  The project
involves a number of experts in formal languages. They currently estimate that the project is about 65\% complete.  As a byproduct
of the FLYSPECK project, Hales, jointly with five coauthors involved in the project, published recently another article  \cite{H10}  on the topic
of the Kepler conjecture, revising the originally published text.

\subsection{Packing space with congruent ellipsoids.}  The problem of packing space with ellipsoids is in sharp contrast with the analogous
two-dimensional problem. In the plane, the density of any packing consisting of congruent ellipses, or even ellipses of equal areas
(see L.~Fejes T\'{o}th \cite{F50}, see also \cite{F72}),  cannot exceed the circle's packing density $\frac{\pi}{\sqrt{12}}$.   It has been noticed in
\cite{BK91} that ellipsoids $E$ exist whose packing density is greater than that of a ball, that is, $\delta(E)>\frac{\pi}{\sqrt{18}}$. The first ellipsoid
found that had this packing property was quite elongated, of a very high {\em aspect ratio}, that is, the ratio of its longest semiaxis to its shortest.

As an improvement of this construction, Wills \cite{W91} found a denser ellipsoid packing, with ellipsoids of a slightly smaller aspect ratio.
However, a much more substantial improvement came about a few years ago.  A.~Donev, F.H.~Stillinger, P. M.~Chaikin,
and S.~Torquato \cite{DSCT04} constructed a remarkably dense packing of congruent ellipsoids that do not differ from a sphere too much,
namely with aspect ratio of $\sqrt3$ (or any greater than that). The packing they found using a computerized experimental simulation technique
reaches density of $0.770732$.  This is the currently highest known density of a packing of space with congruent ellipsoids.

It is not known, however, whether or not there is an upper bound below $1$ for such density.  While no ellipsoid can tile space by its congruent
replicas, thus the packing density of any ellipsoid is smaller than $1$, it is conceivable that an ellipsoid with sufficiently high aspect ratio could
have its packing density as close to $1$ as desired.

\subsection{Packing space with congruent cylinders.} The first non-trivial case of a convex (though unbounded) solid whose packing density,
allowing all isometries, was computed, was the circular cylinder, infinitely long in both directions, that is, the set
$\{(x,y,z)\in{\mathbb R}^3: x^2+y^2\le1\}$, see \cite{BK90a}.  As expected, the maximum density is reached when all cylinders in the packing
are parallel to each other and the plane cross-section of the packing perpendicular to the cylinders forms the densest circle packing in the plane.
In other words, the packing density of the infinite circular cylinder is $\frac{\pi}{\sqrt{12}}$.

The first non-trivial case of a convex compact solid was resolved by A.~Bezdek \cite{B94} who determined the exact value of the packing density
of the rhombic dodecahedron slightly truncated at one of its trihedral vertices.  Although the packing density of Bezdek's example
can be derived from the now proven Kepler conjecture (the truncated rhombic dodecahedron contains the inscribed sphere), Bezdek's proof
was published before the Kepler conjecture was settled and is independent from it.

The packing density of the circular cylinder $\{(x,y,z)\in{\mathbb R}^3: x^2+y^2\le1,\ \ 0\le z\le h\}$ of finite height $h>0$,  conjectured to be
$\frac{\pi}{\sqrt{12}}$ as well, is not known for any value of $h$.  The difficulty of this conjecture is indicated by an example of a certain elliptical
cylinder that admits a packing of density greater than $0.99$ (see \cite{BK91}), while in any arrangement of the congruent copies of it such that all
their generating segments are parallel to each other, the packing's density cannot exceed $\frac{\pi}{\sqrt{12}}$.

Related to the above problem is the following question about tiling space with congruent  right cylinders  (a cylinder is said to be {\em right} if its
generating segment is perpendicular to the plane of its base):

{\em If a right cylinder with a convex base admits a tiling of ${\mathbb R}^3$ with its congruent replicas, must its base admit a tiling of the plane?}

\begin{center}
\begin{figure}[h]
\includegraphics[scale=.8]{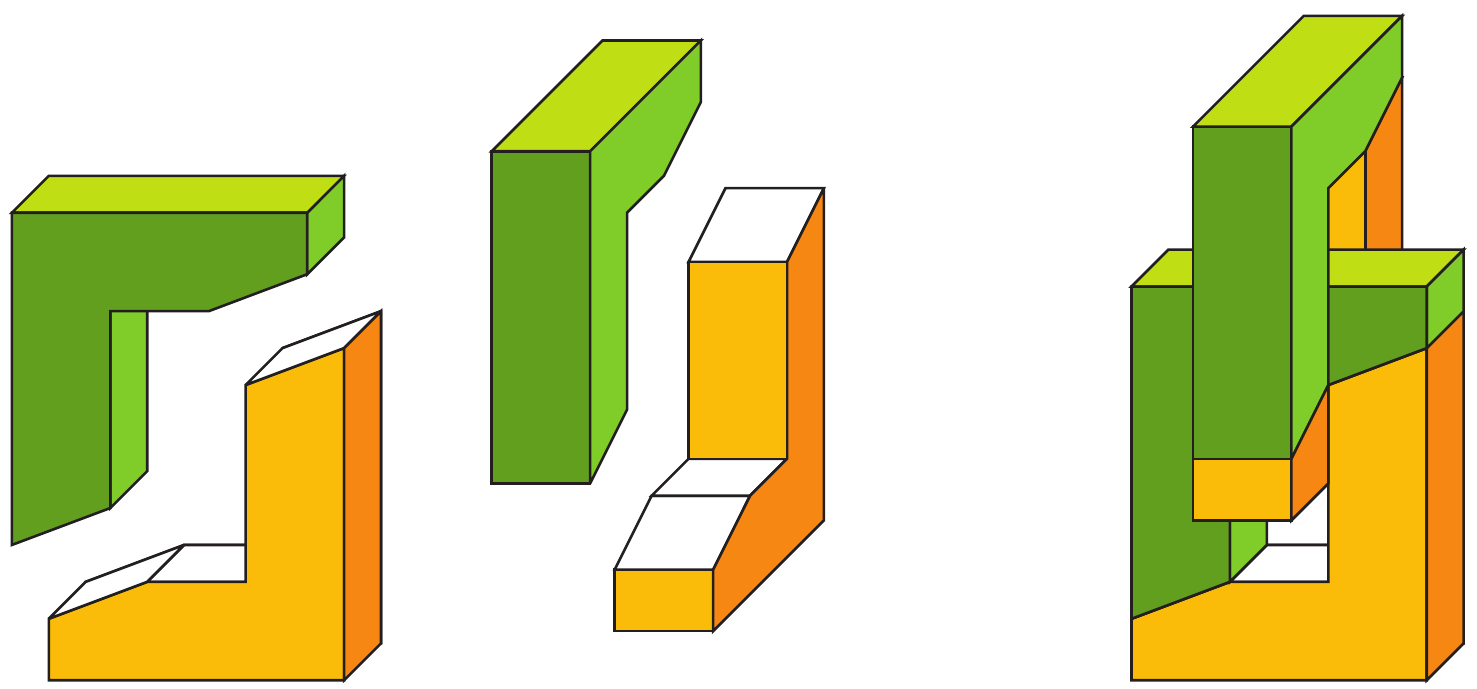}
\caption{A non-convex right prism that tiles space, with base that does not tile the plane.}
\label{cyl1}
\end{figure}
\end{center}

\begin{center}
\begin{figure}[h]
\includegraphics[scale=.8]{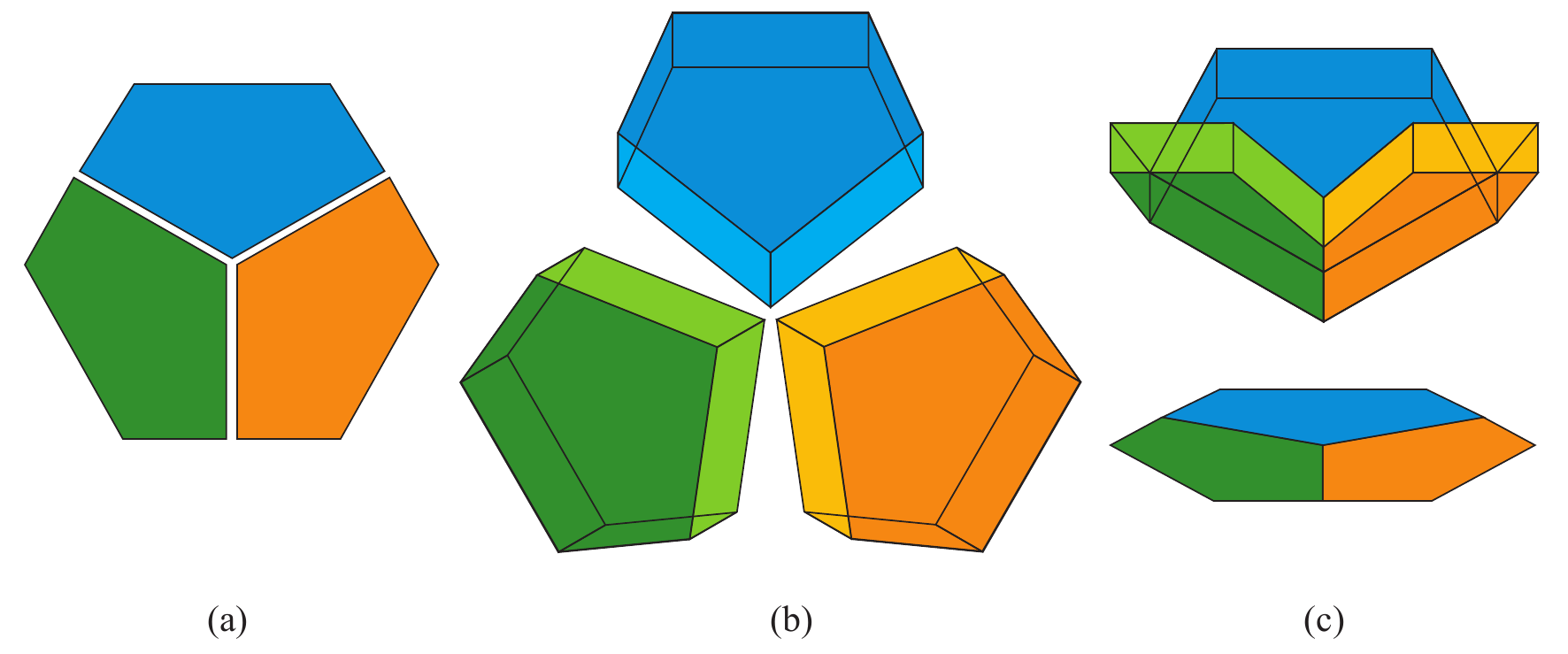}
\caption{A convex, slightly skew prism that tiles space, with base that does not tile the plane.}
\label{cyl2}
\end{figure}
\end{center}

The difficulty of this question is illustrated by two examples from \cite{BK90b}. First, there exists a space-tiling right cylinder with a non-convex
polygonal base that cannot tile the plane (see Fig.~\ref{cyl1}).
Second, there exists a skew prism (though as close to being right as we want) with a convex polygonal base that tiles space,
but whose base cannot tile the plane  (see Fig.~\ref{cyl2}).

To explain the construction shown in Fig.~\ref{cyl2}:  (a)  A regular hexagon is cut into three congruent, axially symmetric pentagons.  (b) By an
affine transormation, the pentagons are stretched slightly, each in the direction of its axis of symmetry, so that they cannot tile the plane.  Then a skew
pyramid is raised over each of the pentagons, so that when they are joined as shown in (c), they form a ``hexagonal cup'' whose projection to the
plane of the original hexagon coincides with the hexagon.  Such ``cups'' can be stacked, forming an infinite beam whose perpendicular cross-section is
the original regular hexagon.  Finally, such hexagonal parallel beams can fill space by the same pattern as the regular hexagon tiles the plane.


The second example shows that the packing density of a cylinder ${\bf\Pi}_s(K)$ over a convex disk $K$ could be greater than the packing density of $K$
in the plane.  In this example, however, the cylinder is skew.  The same phenomenon, however, can occur with a right cylinder as well, as we saw it on the example of a right elliptical cylinder whose packing density is greater than $0.99$.  It seems natural to ask: 

{\em Which convex cylinders have their packing density in space the same as that of their bases in the plane, and which can be packed denser?}

\subsection{High density packing with congruent regular tetrahedra.}  Since no integer multiple of the dihedral angle $\varphi=\arccos\frac{1}{3}=1.23\ldots$ formed by the faces of the regular tetrahedron $\Delta$ equals $2\pi$ ($5\varphi=6.15\ldots$\ \ is just slightly smaller than $2\pi$), we know that $\delta(\Delta)<1$.  Then, how densely can space be packed with congruent regular tetrahedra?  The question is of interest in areas other than mathematics as well, {\em e.g.} physics (compacting loose particles), chemistry (material design), etc.  The past four years brought an exciting development:
a series of articles appeared, each providing a surprisingly dense---denser than previously known---packing.

\begin{center}
\begin{figure}[h]
\includegraphics[scale=.9]{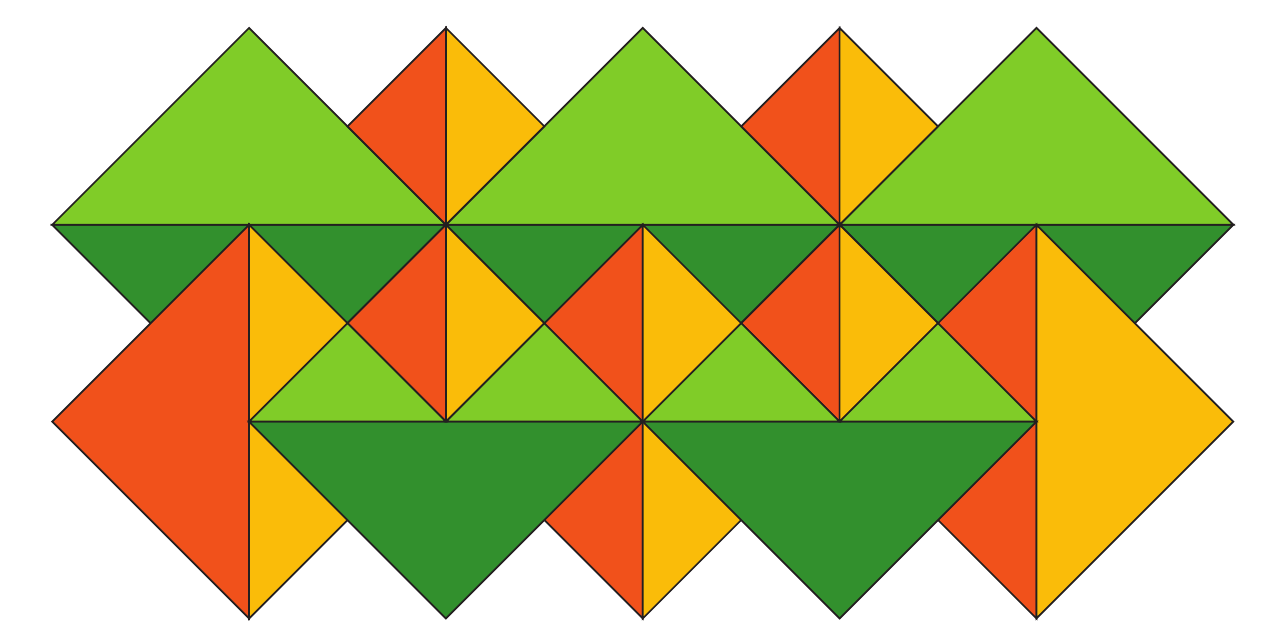}
\caption{A portion of the Conway and Torquato uniform packing of regular tetrahedra. Density: $\frac{2}{3}$.}
\label{CoTo}
\end{figure}
\end{center}

2006. Conway and Torquato \cite{CT06} initiate the race by presenting a surprisingly dense packing with density $0.717455\ldots$, almost twice the lattice packing density of the tetrahedron (see Table 1).  The packing is a lattice arrangement in which the ``repeating unit'' is a cluster of $17$ congruent regular tetrahedra.  Conway and Torquato also give a simple, uniform packing with density $\frac{2}{3}$ (here ``uniform'' means possessing a group of symmetry that acts transitively on the tetrahedra).  This simple packing is a lattice arrangement in which the repeating unit consists of a pair of regular tetrahedra, one rotated by $\frac{\pi}{2}$ with respect to the other (see Fig.~\ref{CoTo}).


Same year, shortly after the appearance of Conway and Torquato's article, Chaikin, Jaoshvili, and Wang  \cite{CWJ07},  a team composed of two
physicists and a high-school student, announce results of an experiment with material tetrahedral dice, packing them tightly, but randomly in
spherical and cylindrical containers.  The experimental results indicate that the packing density of the regular tetrahedron should exceed $0.74$,
perhaps even  $0.76$.

2008.  Elizabeth R.~Chen \cite{C08}, a graduate student at the University of Michigan, Ann Arbor, produces a  packing reaching density $0.7786$,
well above the packing density of the ball.

2009.  Torquato and Jiao \cite{TJ09a, TJ09b}, using computer simulation based on their ``adaptive cell shrinking scheme'' raise Chen's record first to
$0.782\ldots$ and shortly thereafter to $0.823\ldots\,$.

At this point one could hardly expect or predict any significant improvements, but they kept
coming without much delay.

2009.  Haji-Akhbari {\em et al.} \cite{HE09}, using thermodynamic computer simulations that allow a system of particles to evolve naturally towards
high-density states, find a packing whose density reaches $0.8324$.  

2009. Kallus, Elser, and Gravel \cite{KEG10} produce a surprisingly simple uniform one-parameter family of packings - a lattice arrangement of a repeating
unit consisting of just four regular tetrahedra, one pair of tetrahedra joined by a common face and another pair a point-symmetric reflection of the first.
New density record: $\frac{100}{117}=0.85470\ldots\,$.  The packings, though found with the aid of computer, are described analytically.

2010. Torquato and Jiao \cite{TJ10}  produce an analytically described packing with regular tetrahedra bettering the density record of Kallus {\em et al.}  Density: $\frac{12250}{14319}=0.855506\ldots\,$.

2010.  Chen, Engel, and Glotzer \cite{CEG10} set the most recent density record, reached by an analytically described packing.  The currently highest known density is raised to $\frac{4000}{4671} = 0.856347\ldots\,$.

The last few density improvements seem to be inching towards its maximum value.  Though it is difficult to conjecture what that value should be, any 
reasonable upper bound would be welcome as a valuable contribution.  Disappointingly, thus far no specific upper bound, not even by a miniscule
amount  below $1$, has been established.

\end{document}